\documentclass[epj]{svjour}
\RequirePackage{graphicx}
\RequirePackage[ngerman,american]{babel}
\RequirePackage{ntheorem}
\RequirePackage{amsfonts}
\RequirePackage{amsmath}
\RequirePackage{amssymb}
\RequirePackage{mathptmx}
\RequirePackage[scaled=0.9]{helvet}
\RequirePackage{courier}
\RequirePackage{ifthen}
\RequirePackage{color}
\RequirePackage{listings}
\RequirePackage[spaces,obeyspaces]{url}
\RequirePackage{booktabs}
\RequirePackage{verbatim}
\RequirePackage{threeparttable}
\RequirePackage{hyperref}
\RequirePackage{cite}

\newcommand{\bigo}[1]{\mathcal{O}\left(#1\right)}
\newcommand{\E}{\mathrm{e}}

\newcommand{\iid}{i.i.d.}
\newcommand{\npp}{{\sc Npp}}
\newcommand{\ave}[1]{\ensuremath{\mathrm{E}\left[#1\right]}}
\newcommand{\prob}[1]{\ensuremath{\mathrm{P}\left(#1\right)}}

\newcommand{\tabfrac}[2]{%
	\setlength{\fboxrule}{0pt}%
	\fbox{$\frac{#1}{#2}$}%
}

\makeatletter
\def\sgn{\mathop{\operator@font sgn}\nolimits}
\def\srg{\mathop{\operator@font srg}\nolimits}
\def\artanh{\mathop{\operator@font artanh}\nolimits}
\def\EXP{\mathop{\operator@font EXP}\nolimits}
\def\med{\mathop{\operator@font P_{50}}}
\makeatother

\graphicspath{{../figures/}}

\begin{document}

\title{Analysis of the Karmarkar-Karp Differencing Algorithm}

\author{Stefan Boettcher\inst{1} \and Stephan Mertens\inst{2,3}}         
\institute{Department of Physics, Emory University, Atlanta GA 30322-2430, U.S.A. \and Institut f\"ur Theoretische Physik, Otto-von-Guericke Universit\"at, PF 4120, 39016 Magdeburg, Germany \and
Santa Fe Institute, 1399 Hyde Park Road,
Santa Fe, New Mexico 87501, U.S.A. }
\date{\today}
\abstract{
  The Karmarkar-Karp differencing algorithm is the best known
  polynomial time heuristic for the number partitioning problem,
  fundamental in both theoretical computer science and statistical
  physics. We analyze the performance of the differencing algorithm on
  random instances by mapping it to a nonlinear rate equation.  Our
  analysis reveals strong finite size effects that explain why the
  precise asymptotics of the differencing solution is hard to
  establish by simulations. The asymptotic series emerging from the
  rate equation satisfies all known bounds on the Karmarkar-Karp
  algorithm and projects a scaling $n^{-c\ln n}$, where
  $c=1/(2\ln2)=0.7213\ldots$. Our calculations reveal subtle relations
  between the algorithm and Fibonacci-like sequences, and we establish
  an explicit identity to that effect.
 % end of PACS codes
} %end of abstract
\maketitle
\section{Introduction}
\label{sec:intro}

Consider a list of $n$ positive numbers. Replacing the two largest
numbers by their difference yields a new list of $n-1$ numbers.
Iterating this operation $n-1$ times leaves us with a single number.
Intuitively we expect this number to be much smaller than all the
numbers in the original list.  But how small? This is the
question that we address in the present paper.

The operation that replaces two numbers in a list by their difference
is called \emph{differencing}, and the procedure that iteratively
selects the two \emph{largest} numbers for differencing is known as
\emph{largest differencing method} or LDM. This method was introduced
in 1982 by Karmarkar and Karp \cite{karmarkar:karp:82} as an algorithm
for solving the number partitioning problem (\npp): Given a list
$a_1,a_2,\ldots,a_n$ of positive numbers, find a partition, i.e.\ a
subset $A\subset\{1,\ldots,n\}$ such that the discrepancy
\begin{equation}
  \label{eq:cost-function}
  D(A) =
   \Big|\sum_{i\in A}a_i - \sum_{i\not\in A} a_i\Big|,
\end{equation}
is minimized. Obviously, LDM amounts to deciding iteratively that the
two largest numbers will be put on different sides of the partition,
but to defer the decision on what side to put each number. The final
number then represents the discrepancy.

Despite its simple definition, the \npp\ is of considerable
importance both in theoretical computer science and statistical physics.
The \npp\ is NP-hard, which means (a) that no algorithm is known that
is essentially faster than exhaustively searching through all
$2^n$ partitions, and (b) that the \npp\ is computationally
equivalent to many famous problems like the Traveling Salesman
Problem or the Satisfiability Problem \cite{mertens:moore:pc}. In fact,
the \npp\ is one of Garey and Johnson's six basic NP-hard problems
that lie at the heart of the theory of NP-completeness
\cite{garey:johnson:79}, and it is the only one of
these problems that actually deals with numbers.  Hence it is often
chosen as a base for NP-hard\-ness proofs of other problems
involving numbers, like bin packing, multiprocessor scheduling
\cite{bauke:mertens:npp}, quadratic programming or knapsack problems.
The \npp\ was also the base of one of the first public key crypto systems
\cite{merkle:hellman:78}.

In statistical physics, the significance of the \npp\ results from the fact
that it was the first system for which the local REM scenario was
established \cite{mertens:00,borgs:chayes:pittel:01}. 
The notion local REM scenario refers to systems which locally (on the
energy scale) behaves like Derrida's random energy model 
\cite{derrida:80,derrida:81}. It is conjectured to be a universal
feature of random, discrete systems \cite{rem2}.  Recently,
this conjecture has been proven for several spin glass models
\cite{bovier:kurkova:05,bovier:kurkova:07} and for directed polymers
in random media \cite{kurkova:08}.

Considering the NP-hardness of the problem it is no surprise that LDM
(which runs in polynomial time) will generally not find the optimal
solution but an approximation. Our initial question asks for the
quality of the LDM solution to \npp, and to address this question we
will focus on random instances of the \npp\ where the numbers $a_j$
are independent, identically distributed (\iid ) random numbers,
uniformly distributed in the unit interval. Let $L_n$ denote the
output of LDM on such a list. Yakir \cite{yakir:96} proved that the
expectation $\ave{L_n}$ is asymptotically bounded by
\begin{equation}
  \label{eq:yakir}
  n^{-b\ln n} \leq \ave{L_n} \leq n^{-a\ln n}\,,
\end{equation}
where $a$ and $b$ are (unknown) constants such that
\begin{equation}
  \label{eq:yakir2}
  b \geq a \geq \frac{1}{2\ln 2} = 0.7213\ldots\,.
\end{equation}
In this contribution we will argue that $b=a=\frac{1}{2\ln2}$.

The paper is organized as follows. We start with a comprehensive
description of the differencing algorithm, a simple (but wrong)
argument that yields the scaling \eqref{eq:yakir} and a presentation
of simulation data that seems to violate the asymptotic bound
\eqref{eq:yakir2}.  In section \ref{sec:recursions} we reformulate LDM
in terms of a stochastic recursion on parameters of exponential
variates. This recursion will then be simplified to a deterministic, nonlinear
rate equation in section \ref{sec:rateeq}. A numerical investigation of this
rate equation reveals a structure in the dynamics of LDM that can be used as
an Ansatz to simplify both the exact recursions and the rate equation.
This will lead to a simple, Fibonacci like recursion (section \ref{sec:fibonacci}) and to an analytic solution of the rate equation (section \ref{sec:continuum}). In both cases we can derive the asymptotics including the corrections
to scaling, and we claim that a similar asymptotic expansion holds for the original LDM. The latter claim is corroborated by fitting the asymptotic expansion
to the available numerical data on LDM.

% Considered as an optimisation problem, the REM is as hard as it can
% get.  Locating the minimum in a disordered table obviously requires to
% search the whole table, which for the \npp\ grows exponentially with
% $n$.  Heuristic algorithms that yield near optimal solutions by
% considering only a polynomial number of configurations usually rely on
% correlations between configurations and energies, but in a REM no such
% correlations exist. The local REM nature of the \npp\ accounts for the
% poor performance of heuristic algorithms that has been observed in the
% past .

\section{Differencing Algorithm}
\label{sec:differencing}

\begin{figure}[htbp]
  \centering
  \includegraphics[width=0.8\linewidth]{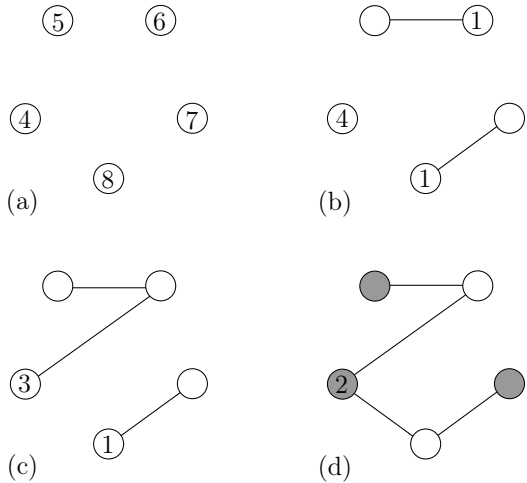}
  \caption{The differencing algorithm in action.}
  \label{fig:differencing}
\end{figure}

The differencing scheme as described in the introduction gives the
value of the discrepancy, but not the actual partition. For that
we need some additional bookkeeping, which is most easily implemented
in terms of graphs (Fig.~\ref{fig:differencing}). The algorithm
maintains a list of rooted trees where each root is
labeled with a number.  The algorithm starts with $n$ trees of size
one and the roots labeled with the numbers $a_i$. Then the following
steps are iterated until a single rooted tree of size $n$ remains:
\begin{enumerate}
\item Among all roots, find those with the largest ($x$) and second largest
($y$) label.
\item Join nodes $x$ and $y$ with an edge, declare node $x$ as the root
  of the new tree and relabel it with $x-y$.
\end{enumerate}
After $n-1$ iterations all nodes are spanned by a tree whose root is labeled
by the final discrepancy. This tree can easily be two-colored, and the colors
represent the desired partition.

Fig.~\ref{fig:differencing} illustrates this procedure on the instance
$(4,5,6,7,8)$. The final two coloring corresponds to the partition
$(4,5,7)$ versus $(6,8)$ with discrepancy $2$. Note that the optimum partition
$(4,5,6)$ versus $(7,8)$ achieves discrepancy $0$.

Technically, LDM boils down to deleting items from and inserting items into a
sorted list of size $n$. This can be done in time $\bigo{n\ln n}$
using an advanced data structure like a heap \cite{kleinberg:tardos}.
Hence LDM is very efficient, but how good is it? As we have already
seen in the example, LDM can miss the optimal partition. And for
random instances, the corridor \eqref{eq:yakir} is far above 
the true optimum, which is known to scale like
$\Theta(\sqrt{n}\,2^{-n})$ \cite{borgs:chayes:pittel:01}. 
Yet LDM yields the best results 
that can be achieved in polynomial time.  Many alternative
algorithms have been investigated in the past
\cite{johnson:etal:91,ruml:etal:96}, but they all produce results
worse than \eqref{eq:yakir}.  The few algorithms that can actually
compete with the Karmarkar-Karp procedure use the same elementary
differencing operation \cite{korf:98,storer:etal:96}. It seems as if
the differencing scheme marks an inherent barrier for polynomial time
algorithms.

The following argument explains the scaling
\eqref{eq:yakir}. The typical distance between adjacent pairs of the
$n$ numbers in the interval $[0,1]$ is $n^{-1}$. Hence after $n/2$
differencing operations we are left with $n/2$ numbers in the interval
$[0,n^{-1}]$. The typical distance between pairs is now $2n^{-2}$.
After another round of $n/4$ differencing operations we get $n/4$
numbers in the range $[0,8n^{-3}]$.  In general, after $2^k$
differencing operations we are left with $n/2^k$ numbers in the range
$[0,2^{k \choose 2}n^{-k}]$.  Reducing the original list to a single
number requires $k=\log_2 n$ differencing operations, and applying the
above argument all the way down suggests that
\begin{equation}
  \label{eq:simple-argument}
  \ave{L_n} \propto n^{-c\ln n}
\end{equation}
with
\begin{equation}
  \label{eq:simple-c}
  c = \frac{1}{2 \ln 2} = 0.721\ldots\,.
\end{equation}
As we will see, this is the right scaling, yet the argument cannot be
correct. This follows from the fact that it predicts the same scaling
for the paired differencing method (PDM). Here in each round all pairs
of adjacent numbers are replaced by their difference in parallel.
This method, however, yields an average discrepancy of order
$\Theta(n^{-1})$ \cite{lueker:87}. Yet, our analysis below suggests
that \eqref{eq:simple-argument} and \eqref{eq:simple-c} indeed
describe the asymptotic behavior correctly, although a far more
subtle treatment is required.

% \section{Monte Carlo Studies of LDM}
% \label{sec:MC}

\begin{figure}
  \centering
  \includegraphics[angle=270,width=\columnwidth]{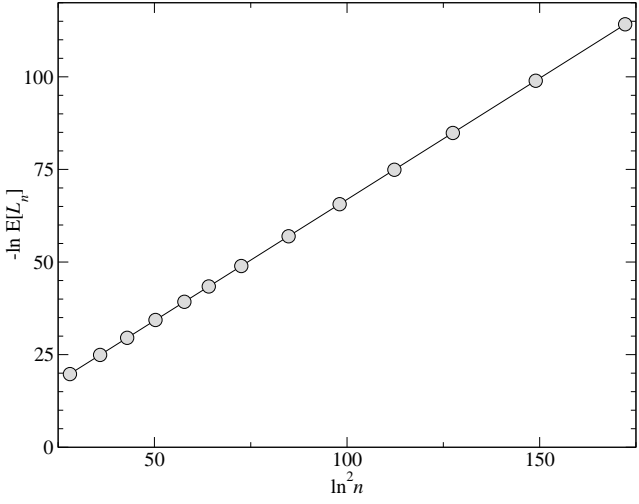}
  \caption{Results of LDM applied to $n$ random \iid\ numbers, uniformly
drawn from the unit interval. Each data point represents between $10^5$ (large $n$) and $10^7$ samples (small $n$). The solid line is the linear fit
$-\ln{\ave{L_n}} = 1.42 + 0.65\,\ln^2n$.}
  \label{fig:KKsim}
\end{figure}

An obvious approach to find the quality of LDM are simulations.  We
ran LDM on random instances of varying size $n$, and Figure
\ref{fig:KKsim} shows the results for $\ave{L_n}$.  Apparently
$\ln\ave{L_n}$ scales like $\ln^2 n$, in agreement with \eqref{eq:yakir}
and \eqref{eq:simple-argument}. A linear fit seems to yield
\begin{displaymath}
  c \simeq 0.65 
\end{displaymath}
for the constant in \eqref{eq:simple-argument}, which clearly violates
the bound $c \geq 1/2\ln 2$. Apparently even $n=10^6$ is too small
to see the true asymptotic behavior. This may be the reason why
Monte Carlo studies of LDM never have been published.

\begin{figure}
  \centering
  \includegraphics[angle=270,width=\columnwidth]{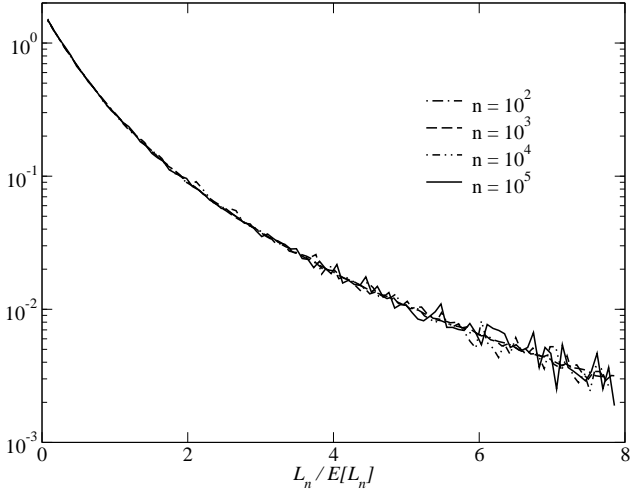}
  \caption{Probability density function of $L_n/\ave{L_n}$.}
  \label{fig:KKpdf}
\end{figure}

A plot of the probability density function (pdf) of $L_n/\ave{L_n}$
reveals a data collapse varying values of $n$ (Fig.~\ref{fig:KKpdf}).
Apparently the complete statistics of $L_n$ is asymptotically
dominated by a single scale $n^{-c\ln n}$.

Some technical notes about simulating LDM are appropriate.
Differencing means subtracting numbers over and over again. The
numerical precision must be adjusted carefully to support this and to
be able to represent the final discrepancy of order $n^{-c\ln n}$. We
used the freely available GMP library \cite{gmp} for the required
multiple precision arithmetic and ran all simulations on $\ell$-bit
integers where the number of bits ranges from $\ell=40$ (for $n=20$)
to $\ell=300$ for $n=1.5\cdot 10^7$. The integer discrepancies were then
rescaled by $2^{-\ell}$. The pseudo random number generator was taken
from the TRNG library \cite{trng}.

\section{Exact Recursions}
\label{sec:recursions}

A common problem in the average-case analysis of algorithms like LDM
is that numbers become conditioned and cease to be independent as 
the algorithm proceeds. Lueker \cite{lueker:87} proposed to use
exponential instead of uniform variates to cope with this problem.
Let $X_1,\ldots,X_{n+1}$ be \iid\ random exponentials with mean $1$
and consider the partial sums $S_k = \sum_{i=1}^k X_i$. Then the joint
distribution of the ratios $S_k/S_{n+1}$, $k=1,\ldots,n$, is the same
as that of the order statistics of $n$ \iid\ uniform variates from
$[0,1]$ \cite{feller:vol2}.
% section III.3, examples (d) and (e)
As a consequence, LDM will produce
the same distribution of data no matter whether it is run on uniform variates
or on $S_k/S_{n+1}$. Let $\hat{L}_n$ denote the
result of LDM on the partial sums $S_1, S_2, \ldots, S_n$.
Since the output of LDM is linear in its input, we have
\begin{equation}
  \label{eq:L-hat}
  \hat{L}_n \stackrel{D}{=} S_{n+1} L_n\,,
\end{equation}
where $S_{n+1}$ is the sum of $n+1$ \iid\ exponential variates and the
notation $X\stackrel{D}{=}Y$ indicates that the random variable $X$
and $Y$ have the same distribution. The probability density of
$S_{n+1}$ is the gamma density
\begin{equation}
  \label{eq:gamma-density}
  g_{n+1} (s) = \frac{s^n}{n!}\,\E^{-s}\,.
\end{equation}
Taking expectations of both sides of \eqref{eq:L-hat}
we get
\begin{equation}
  \label{eq:exp-relation}
  \ave{L_n} = \frac{\ave{\hat{L}_n}}{n+1}\,.
\end{equation}
This allows us to derive the asymptotics of
$\ave{L_n}$ from the asymptotics
of $\ave{\hat{L}_n}$.

Exponential variates are well suited for the analysis of LDM because
the sum and difference of two exponential variates are again
exponential variates. Once started on exponential variates, LDM
keeps working on exponentials all the time. This allows us to express
the operation of LDM in terms of a recursive equation for the
parameters of exponential densities \cite{yakir:96}. 
We start with the following Lemma:
\begin{lemma}
 \label{lem:lem1}
 Let $X_1$ and $X_2$ be independent exponential random variables with
 parameter $\lambda_1$ and $\lambda_2$, resp.. The probability 
 of the event $X_1<X_2$ is given by
 \begin{equation}
   \label{eq:lem1}
   \prob{X_1<X_2} = \frac{\lambda_1}{\lambda_1 + \lambda_2}\,.
 \end{equation}
 Furthermore, conditioned on the event $X_1<X_2$, the variables $X_1$ and $X_2-X_1$ 
 are independent exponentials with parameters $\lambda_1+\lambda_2$
 (for $X_1$) and $\lambda_2$ for $X_2-X_1$.
\end{lemma} 
The proof of Lemma \ref{lem:lem1} consists of trivial integrations of
the exponential densities and is omitted here.

Next we consider generalized partial sums
of exponentials, described by $n$-tuples 
\begin{displaymath}
  (\lambda_1,\lambda_2,\ldots,\lambda_n)\,.
\end{displaymath}
This $n$-tuple is shorthand for the sequence of partial
sums 
\begin{displaymath}
  (X_1,X_1+X_2,\ldots,\sum_{i=1}^n X_i)
\end{displaymath}
with
$X_i=\EXP(\lambda_i)$.

Now let us look at the result of one iteration of LDM on
$(\lambda_1,\lambda_2,\ldots,\lambda_n)$. The two largest numbers
are removed and replaced by their difference $X_n$ which is an
$\EXP(\lambda_n)$ variate. Lemma \ref{lem:lem1} 
tells us, that the probability that this number is
the smallest in the list is
\begin{displaymath}
  \prob{X_n < X_1} = \frac{\lambda_n}{\lambda_1+\lambda_n},
\end{displaymath}
and conditioned on that event, the smallest number is 
an $\EXP(\lambda_1+\lambda_n)$ variate and the increment to the
2nd smallest number $X_1-X_n$ is an independent
$\EXP(\lambda_1)$ variate. Conditioned on $X_n < X_1$ we get
another $\lambda$-tuple as the input for the next iteration: 
\begin{displaymath}
  X_n < X_1 \Rightarrow (\lambda_1+\lambda_n, \lambda_1, \lambda_2, \ldots, \lambda_{n-2})
\end{displaymath}
The probability that $X_n \geq X_1$ is
\begin{displaymath}
  \prob{X_n \geq X_1} = \frac{\lambda_1}{\lambda_1+\lambda_n},
\end{displaymath}
and in this case $X_1$ is an $\EXP(\lambda_n+\lambda_1)$ variate,
whereas the difference $X_n - X_1$ is an $\EXP(\lambda_n)$ variate.
Now the probability that the new number $X_n$ is second in the
new list reads
\begin{eqnarray*}
  \prob{X_n \geq X_1 \cap X_n < X_1+X_2} &=& \prob{X_n \geq X_1 \cap X_n-X_1 < X_2}\\
                                 &=& \frac{\lambda_1}{\lambda_1+\lambda_n}
                                     \frac{\lambda_n}{\lambda_2+\lambda_n}
\end{eqnarray*}
and conditioned on that event the input for the new iteration is
\begin{displaymath}
  (\lambda_1+\lambda_n, \lambda_2 + \lambda_n, \lambda_2, \ldots, \lambda_{n-2})\,.
\end{displaymath}
This argument can be iterated to calculate the probability of
$X_n$ becoming the $k$-th number in the new list. Denoting the
partial sums by $S_k$ we get
\begin{equation}
  \label{eq:P-k1}
  \prob{X_n \geq S_{k-1} \cap X_n < S_k} = \frac{\lambda_n}{\lambda_k+\lambda_n}
  \prod_{i=1}^{k-1} \frac{\lambda_i}{\lambda_i+\lambda_n}
\end{equation}
for $k=1,\ldots,n-2$ and conditioned on that event the new list is
\begin{equation}
  \label{eq:P-lambda1}
  (\lambda_1+\lambda_n, \ldots, \lambda_k+\lambda_n, \lambda_k, \lambda_{k+1}, \ldots, \lambda_{n-2})\,.
\end{equation}
The final case is that $X_n$ becomes the largest number in the new list.
This happens with probability
\begin{equation}
  \label{eq:P-k2}
  \prob{X_n \geq S_{n-2}} = \prod_{i=1}^{n-2} \frac{\lambda_i}{\lambda_i+\lambda_n}
\end{equation}
and leads to the list
\begin{equation}
  \label{eq:P-lambda2}
  (\lambda_1+\lambda_n, \ldots, \lambda_{n-2}+\lambda_n, \lambda_n)\,.
\end{equation}
In all cases we stay within the set of instances given by 
partial sums of independent exponentials, and we can apply 
Eqs.~\eqref{eq:P-k1} to \eqref{eq:P-lambda2} recursively until
we have reduced the original problem to a 
$(\lambda_1, \lambda_2)$-instance which tells us that the final
difference is an $\EXP(\lambda_2)$ variate.

 \begin{figure}[htbp]
   \centering
   \includegraphics[width=0.8\columnwidth]{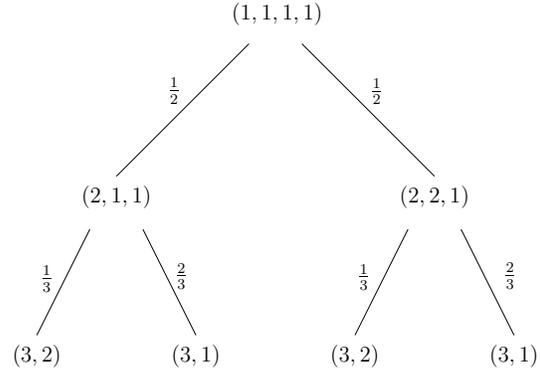}
   \caption{Statistics of LDM on $n=4$. The final difference is distributed
 according to $p_4(x)=\frac{2}{3}\cdot \E^{-x} + \frac{1}{3}\cdot 2\E^{-2x}$}
   \label{fig:kktree}
 \end{figure}

\begin{table}[htb]
  \centering
  \begin{tabular}{r|ccccc}
    $k \backslash n$ & 4 & 5 & 6 & 7 & 8 \\\hline
     1 & $\tabfrac{2}{3}$ & $\tabfrac{13}{24}$ & $\tabfrac{41}{120}$ & $\tabfrac{49}{180}$ & $\tabfrac{431}{2520}$ \\
     2 & $\tabfrac{1}{3}$ & $\tabfrac{1}{6}$ & $\tabfrac{5}{18}$ & $\tabfrac{1}{8}$ & $\tabfrac{527}{3456}$ \\
     3 &  & $\tabfrac{7}{24}$ & $\tabfrac{7}{72}$ & $\tabfrac{1073}{4320}$ & $\tabfrac{3079}{38880}$ \\
     4 & & & $\tabfrac{41}{180}$ & $\tabfrac{47}{720}$ & $\tabfrac{1229}{5600}$ \\
     5 & & & $\tabfrac{1}{18}$ & $\tabfrac{53}{360}$ & $\tabfrac{149}{2100}$ \\
     6 & & & & $\tabfrac{7}{72}$ & $\tabfrac{486359}{5443200}$ \\
     7 & & & & $\tabfrac{161}{4320}$ & $\tabfrac{343}{4320}$ \\
     8 & & & & $\tabfrac{1}{135}$ & $\tabfrac{11}{144}$ \\
     9 & & & & & $\tabfrac{26083}{604800}$ \\
    10 & & & & & $\tabfrac{859}{77760}$ \\
    11 & & & & & $\tabfrac{941}{155520}$ \\
    12 & & & & & $\tabfrac{1}{1050}$ \\
    13 & & & & & $\tabfrac{1}{1800}$
  \end{tabular}
  \caption{Coefficients $a_k^{(n)}$ in \eqref{eq:pn}.}
  \label{tab:akn}
\end{table}

Fig.~\ref{fig:kktree} shows the result of this analysis on the
input $(1,1,1,1)$, our original problem with $n=4$. We have
to explore the tree that branches according to the position that
is taken by the new number inserted in the shortened list. The
numbers written on the edges of the tree are the probabilities
for the corresponding transition. Note that we have combined the two
branches emerging from the root that both lead to a $(2,2,1)$-configuration
into a single one by adding their probabilities.
In the end we get
\begin{displaymath}
  p_4(x) = \frac{2}{3}\E^{-x} + \frac{2}{3}\E^{-2x}
\end{displaymath}
for the probability density function (pdf) of $\hat{L}_4$.
In general, the pdf of $\hat{L}_n$ is a sum of exponentials,
\begin{equation}
  \label{eq:pn}
  p_n(x) = \sum_{k} a_k^{(n)} k\, \E^{-kx}
\end{equation}
where $a_k^{(n)}$ is the probability of LDM returning an
$\EXP(k)$-variate.  For small values of $n$, this probabilities
can be calculated by expanding the recursions explicitly
(Table \ref{tab:akn}), but for larger values of $n$ this approach
is prohibited by the exponential growths of the number $K(n)$ of branches
that have to be explored. 
% $K(n)$ is the largest value of $k$ such that
% $a_k^{(n)} > 0$. Then $K(n)$ is a Fibonacci sequence,
% \begin{equation}
%   \label{eq:fibo-kn}
%   K(n) = K(n-1) + K(n-2)
% \end{equation}
% with $K(3)=K(2)=1$. This implies that
% $K(n) \simeq \Phi^n$, where $\Phi=1.618\ldots$ is the
% golden ratio.

\begin{figure}
  \centering
  \includegraphics[angle=270,width=\columnwidth]{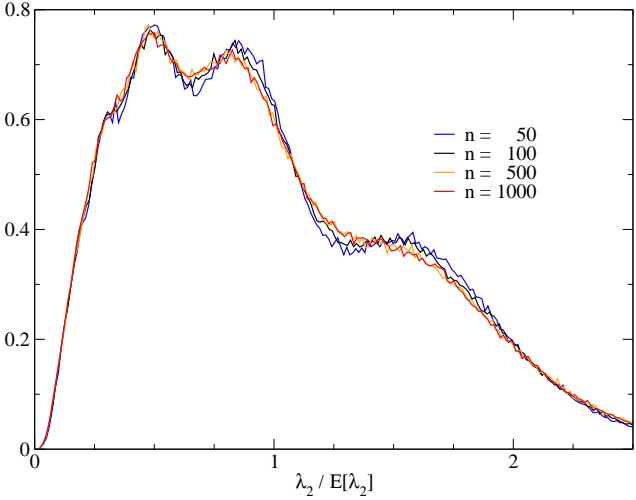}
  \caption{Probability density function of $\lambda_2/\ave{\lambda_2}$.}
  \label{fig:kk-walk-pdf}
\end{figure}

Alternatively we can explore the tree of $\lambda$-tuples by walking
it randomly.  Given a tuple $(\lambda_1\ldots,\lambda_n)$, we generate
a random integer $1 \leq k \leq n-1$ with probability
\begin{equation}
  \label{eq:draw-k}
  \prob{k \leq \ell} = \begin{cases}
     1 - \prod_{j=1}^\ell \frac{\lambda_j}{\lambda_j+\lambda_n} & (\ell < n-1) \\
     1 & (\ell = n-1)
  \end{cases}
\end{equation}
and using this random $k$ we generate a new tuple of size $n-1$
according to Eqs.~\eqref{eq:P-lambda1} or \eqref{eq:P-lambda2}.
This process is iterated until the tuple size is two, and the
final value of $\lambda_2$ is the parameter for
the statistics of $\hat{L}$. The probability density of 
$\lambda_2/\ave{\lambda_2}$
is shown in Fig.~\ref{fig:kk-walk-pdf}). Again the data collapse 
corroborates the claim that the statistics of $LDM$ is
dominated by a single scale.

\section{Rate Equation}
\label{sec:rateeq}

We can turn the exact recursions from Sec.~\ref{sec:recursions} into a
set of rate equations for the time-evolution of the \emph{average} $\lambda$-tuple. Let $\lambda_i^t$ denote the value of $\lambda_i$ after $t$ iterations, such that
\begin{equation}
\left(\lambda_1^t,\lambda_2^t,\ldots,\lambda_{n-t}^t\right)\to\left(\lambda_1^{t+1},\lambda_2^{t+1},\ldots,\lambda_{n-t-1}^{t+1}\right)\,.
\label{eq:lambda_trans}
\end{equation}
As explained in Sec.~\ref{sec:recursions}, at ``time''
$t$ a number $k$, $1\leq k\leq n-1-t$ is chosen with probability
\begin{equation}
  \label{eq:draw-k-t}
  \mathrm{P}_t\left(k \leq \ell\right) = \begin{cases}
     1 - \prod_{j=1}^\ell \frac{\lambda_j^t}{\lambda_j^t+\lambda_{n-t}^t} & (\ell < n-1-t) \\
     1 & (\ell = n-1-t)
  \end{cases}\,.
\end{equation}

% \begin{eqnarray}
% p_k^t & = & \begin{cases}
% \left(1-\alpha_k^t\right)\prod_{j=1}^{k-1}\alpha_j^t, & 1\leq k<n-1,\\ \\
% \hfill\prod_{j=1}^{n-2}\alpha_j^t, & k=n-1,
% \end{cases}
% \label{eq:kprob}
% \end{eqnarray}
% with $\prod_{j=a}^bf_j\equiv1$ for $b<a$ and  abbreviating
% \begin{equation}
% \alpha_i^t=\frac{\lambda_i^t}{\lambda_i^t+\lambda_{n-t}^t}\,.
% \label{eq:alphadef}
% \end{equation}
% Note that $0\leq\alpha_i^t\leq1$ for all $(i,t)$.

% Given a chosen $k$ at time $t$, the update results in
% \begin{eqnarray}
% \label{eq:lambdacase}
%  & \left(\lambda_1^{t+1},\lambda_2^{t+1},\ldots,\lambda_{n-t-1}^{t+1}\right)=\\
%  & \begin{cases}
% \left(\lambda_1^t+\lambda_{n-t}^t,\ldots,\lambda_k^t+\lambda_{n-t}^t,\lambda_k^t,\ldots,\lambda_{n-t-2}^t\right),
%  & 1\leq k<n-t-1,\\ \\
% \left(\lambda_1^t+\lambda_{n-t}^t,\ldots,\lambda_{n-t-2}^t+\lambda_{n-t}^t,\lambda_{n-t}^t\right), & k=n-t-1.\end{cases}
% \nonumber
% \end{eqnarray}
Depending on the choice of $k$, Eqs.~\eqref{eq:P-lambda1} and
\eqref{eq:P-lambda2} suggest that $\lambda_i^{t+1}$ only takes on one
of \emph{two} possible values. For $1\leq i<n-t-1$, these are
\begin{equation}
\lambda_{i}^{t+1} =  
\begin{cases}
\lambda_{i}^{t}+\lambda_{n-t}^{t} & (i\leq k\leq n-t-1)\\
\lambda_{i-1}^{t} & (1\leq k<i)
\end{cases}\,,
\label{eq:lambdachoice}
\end{equation}
whereas for $i=n-t-1$, the two values are 
\begin{equation}
\label{eq:lambdachoice2}
\lambda_{n-t-1}^{t+1} =  
\begin{cases}
\lambda_{n-t}^{t} & (k=n-t-1) \\
\lambda_{n-t-2}^{t} & (1\leq k<n-t-1)
\end{cases}\,.
\end{equation}
% Using Eq.~(\ref{eq:kprob}), we obtain the probability of $k<i$ by
% \begin{eqnarray}
% \sum_{l=1}^{i-1}p_{l}^{t} & = & \sum_{l=1}^{i-1}\left(1-\alpha_{l}^{t}\right)
% \prod_{j=1}^{l-1}\alpha_{j}^{t}=1-\prod_{j=1}^{i-1}\alpha_{j}^{t}\quad,
% \label{eq:probksmaller}
% \end{eqnarray}
% which  is easily shown by induction.
We introduce the shorthand
\begin{equation}
{\cal P}_{i}^t=\prod_{j=1}^{i-1}\frac{\lambda^t_j}{\lambda^t_j+\lambda^t_{n-t}}\,,
\label{eq:defprob} 
\end{equation}
for the probability of $k\geq i$ at iteration $t$. 
On average, the evolution of $\lambda_i^t$ is given by the
\emph{rate equation}
\begin{equation}
  \label{eq:lambda1eq}
  \lambda_{i}^{t+1} = \lambda_{i-1}^{t}\left(1-{\cal P}_{i}^t\right)
+\left(\lambda_{i}^{t}+\lambda_{n-t}^{t}\right){\cal P}_{i}^t\,,
\end{equation}
for all $1\leq i<n-1-t$, and at the upper boundary
\begin{equation}
  \label{eq:lambda2eq}
\lambda_{n-(t+1)}^{t+1} = 
\lambda_{n-2-t}^{t}\left(1-{\cal P}_{n-1-t}^t\right)
+\lambda_{n-t}^{t}{\cal P}_{n-1-t}^t\,.
\end{equation}
These equations are defined on the triangular domain $0\leq t\leq
n-1$, $1\leq i\leq n-t$.  The initial conditions are
\begin{equation}
\lambda_{i}^{t=0}=1\qquad\left(1\leq i\leq n\right)\,.
\label{eq:lambdainiteq}
\end{equation}
As described in Sec.~\ref{sec:recursions}, the process terminates at
$t=n-2$ with $\lambda_2^{n-2}$ characterizing the exponential variate
for the final difference in LDM. Yet, Eq.~(\ref{eq:lambda2eq}) for
$t=n-2$ implies $\lambda_1^{n-1}=\lambda_2^{n-2}$, reflecting
the final, trivial differencing step, and it will prove
conceptually advantageous to focus on the asymptotic properties of
$\lambda_1^{n-1}$ instead.

% \begin{figure}[htbp]
% \centering
% \includegraphics[angle=270,width=0.95\columnwidth]{rate-eq}
% \caption{The simulations from Fig.~\ref{fig:KKsim}
%   (squares) and the output of the rate
%   equations~(\ref{eq:lambda1eq}-\ref{eq:lambdainiteq}) (circles) approach
%   each other asymptotically. The lines are quadratic fits.
% }
% \label{fig:rate-eq}
% \end{figure}

Since the rate equation is an approximation to the exact
recursion, we need to check how accurate it is. 
We have solved the rate
equations~(\ref{eq:lambda1eq}-\ref{eq:lambdainiteq}) numerically up to
$n =5\cdot 10^6$. 
Fig.~\ref{fig:allextrapolations} shows
\begin{displaymath}
  \frac{\ln\big(\lambda_1^{n-1}\,(n+1)\big)}{\ln^2 n}
\end{displaymath}
from the rate equation versus $1/\ln n$. If $\lambda_1^{n-1}$ were calculated
as an average from the exact recursion, it should be equal to
\begin{displaymath}
  -\frac{\ln\ave{L_n}}{\ln^2n}
\end{displaymath}
from the direct simulation of LDM. Fig.~\ref{fig:allextrapolations} shows this
quantity, too. Apparently the error introduced by approximating the
exact recursion by the rate equation vanishes for $n\to\infty$, and our
conjecture is that the rate equation and the exact recursion 
are asymptotically equivalent. Judging from our numerical studies
below, see Tab.~\ref{tab:fitparameter}, both asymptotic series have a relative
difference of size $\ln\ln n/\ln^2(n)$.

The time to solve the rate equation numerically scales like
$\bigo{n^2}$, so it is actually more efficient to simulate LDM
directly, not least because the sampling for the latter can be done
efficiently on a parallel machine. For analytic approaches, however,
the rate equation is more convenient.

\begin{figure}[htbp]
\centering
\includegraphics[angle=0,width=1.0\columnwidth]{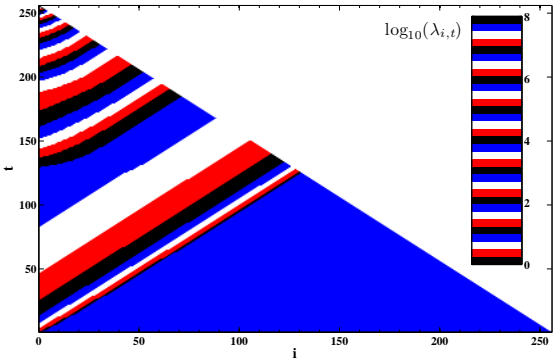}
\caption{Contour plot on a logarithmic scale for the numerical
  solution $\lambda_i^t$ of the rate
equations~(\ref{eq:lambda1eq}-\ref{eq:lambdainiteq}) at $n=256$. The
solution is $\lambda_i^t\simeq 1$ throughout the entire lower triangle, and
it increases monotonically for increasing $t$ above that. The solution
rises by about a decade between each repeat of a band color. Note
the ever more-rapid alternation between narrowing and widening bands, signifying
regions of rapid gain interrupted by extended plateaus. The
regular banded structure along diagonals $t-i=const$ justifies the
similarity solution in Eq.~(\ref{eq:continuumwaveAnsatz}). The only
notable exceptions occur in asymptotically diminishing regions near
$i=1$ and $t=n/2,3n/4,7n/8,\ldots$.
}
\label{fig:lambdacontour}
\end{figure}

\noindent
The initial probabilities decay exponentially,
\begin{equation}
  \label{eq:initial_P}
  \mathcal{P}_i^0 = 2^{-i}\,,
\end{equation}
which implies that only the first values $\lambda_1,\lambda_2,\ldots$
increase. Everywhere else, $\mathcal{P}_i$ is essentially zero, and
those entries will not increase until the first term of
\eqref{eq:lambda1eq} has copied the values from the low index
boundary. Hence we expect a ``wavefront'' of increased $\lambda$-values 
to travel with a velocity one index per time step toward the upper boundary,
which in turn travels with the same velocity towards the lower boundary.
As can be seen from Fig.~\ref{fig:lambdacontour}, this traveling
wavefronts of increasing heights are a hallmark of the rate equation for
all times $t$. We will use this intuitive picture for an Ansatz to analyze
both the exact recursion and the rate equation in the next two sections.

% \begin{figure}[htbp]
% \centering
% \includegraphics[angle=0,width=0.95\columnwidth]{all_extrapolation}
% \caption{Extrapolation plot of various predictions for $\ave{L_n}$ in
%   Eq.~(\ref{eq:cost-function}). The simulations from Fig.~\ref{fig:KKsim}
%   (circles) and the output of the rate
%   equations~(\ref{eq:lambda1eq}-\ref{eq:lambdainiteq}) (squares) approach
%   each other asymptotically. Both extend well (dashed lines) to the
%   scaling in Eq.~(\ref{eq:simple-argument}),
%   $-\ln\left(\ave{L_n}\right)/\ln\left(n^2\right)\sim1/2\ln2=0.7213$
%   (dotted line). Shown are also the approximate results of the
%   sequence in Eq.~(\ref{eq:lambdafibo}) (diamonds) and the asymptotic
%   expansion from Eq.~(\ref{eq:f1asympt}) (full line).
% }
% \label{fig:KKextra}
% \end{figure}

\section{Fibonacci Model}
\label{sec:fibonacci}

Both the exact recursion and the rate equations yield
\begin{equation}
\lambda_{1}^{t+1} = \lambda_{1}^{t}+\lambda_{n-t}^{t}
\label{eq:lambda1recur}
\end{equation}
for the lower boundary that we are ultimately interested in.
This recursion connects the lower and the upper
boundaries at $i=1$ and at $i=n-t$.
Unfortunately, $\lambda_{n-1}^t$ depends in a complicated way on 
entries of the $\lambda$-tuple at different times and different 
places. However, Fig.~\ref{fig:lambdacontour} suggests a 
\emph{similarity Ansatz}
\begin{equation}
  \label{eq:similarity-discrete}
  \lambda_i^t = \lambda_{i-x}^{t-x}\,,
\end{equation}
which makes the upper boundary readily available:
% \begin{figure}
% \includegraphics[scale=0.5]{triangular}
% \caption{Demonstration of Eq.~(\ref{eq:lambda1recur}).
% \label{fig:triangular}}
% \end{figure}
% Fig.~\ref{fig:triangular} demonstrates the
% meaning of Eq.~(\ref{eq:lambda1recur}) for $n=8$. At $t=0$,
% $\lambda_{1}^{1}=\lambda_{1}^{0}+\lambda_{8}^{0}=2$ due to the initial
% conditions $\lambda_{i}^{0}=1$ for all $1\leq i\leq n$.  Then, for
% $t=1$ it is $\lambda_{1}^{2}=\lambda_{1}^{1}+\lambda_{7}^{1}$, but due to
% similarity, all values along each thick diagonal line in
% Fig.~\ref{fig:triangular} are equal, thus
% $\lambda_{7}^{1}=\lambda_{6}^{0}=1$.  In this manner, the recursion in
% Eq.~(\ref{eq:lambda1recur}) picks up (every second of) the initial
% conditions, until $t\geq\left\lfloor n/2\right\rfloor $.  Beyond this
% time, $\lambda_{n-t}^{t}=\lambda_{1}^{2t-n+1}$, for instance,
% $\lambda_{1}^{5}=\lambda_{1}^{4}+\lambda_{4}^{4}=\lambda_{1}^{4}+\lambda_{1}^{1}$
% in Fig.~\ref{fig:triangular}. In general, we find
% \begin{eqnarray}
% \lambda_{1}^{t+1} & = & \lambda_{1}^{t}+\begin{cases}
% \lambda_{n-2t}^{0}, & t<\left\lfloor \frac{n}{2}\right\rfloor ,\\ \\
% \lambda_{1}^{2t-n+1}, & t\geq\left\lfloor \frac{n}{2}\right\rfloor .
% \end{cases}
% \label{eq:lambdafibo}
% \end{eqnarray}
\begin{equation}
  \label{eq:lambdafibo}
  \begin{aligned}
  \lambda_1^{t+1} &= \lambda_1^t + \lambda_1^{2t-n+1} & (0 \leq t < n-1)\\
  \lambda_1^{t} &= 1 & (t \leq 0)
  \end{aligned}
\end{equation}
Note that we have extended the initial conditions
$\lambda_i^t=1$ to hold for all negative times, too.

\begin{figure}
\includegraphics[width=\columnwidth]{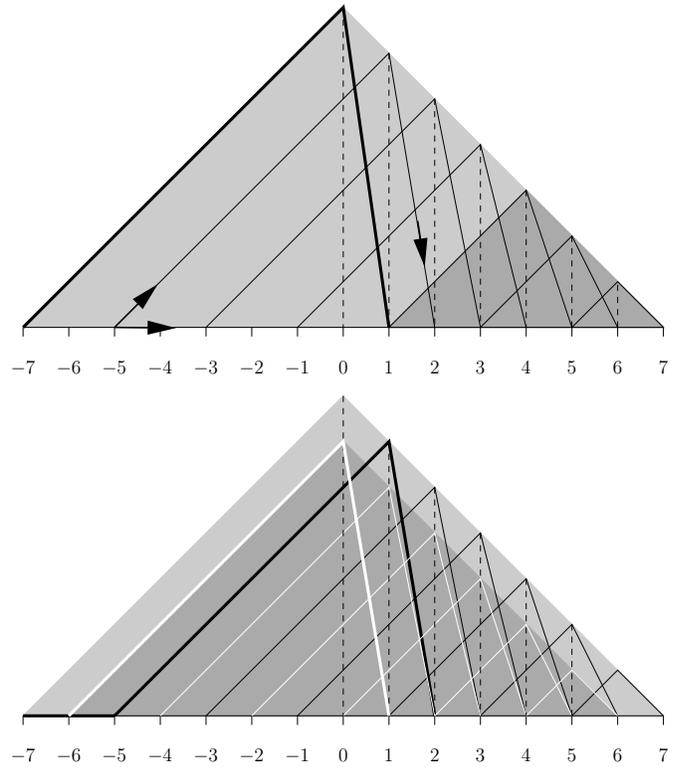}
\caption{Proof of the Fibonacci recursion: The number of different paths
from the leftmost point to the rightmost point in the triangle for $n$ is
the sum of the number of paths in the corresponding triangle of size $[n/2]$
(top) plus the number of paths in the triangle of size $n-1$ (bottom). 
\label{fig:proof}}
\end{figure}

It turns out that one can express the final value $\lambda_1^{n-1}$ of
this recursion in terms of the corresponding values in smaller
systems, which leads to a simple recursion in $n$. To derive this
recursion it is convenient to visualize \eqref{eq:lambdafibo} in terms
of paths in a right-angled triangle $\Delta_n$ (Fig.~\ref{fig:proof}).
The hypotenuse of $\Delta_n$ represents $t$ and ranges from $-n+1$ to
$n-1$, the height is $n-1$. Let us discuss the basic mechanism for the
example $n=8$. The final recursion reads
\begin{displaymath}
  \lambda_1^{7} = \lambda_1^6 + \lambda_1^5\,,
\end{displaymath}
and the two terms on the right hand side correspond to two paths: one
that connects 6 with 7 along the hypotenuse, the other connects 5 with
7 along the path that is ``reflected'' at the right leg of $\Delta_8$.
In our case a reflected path moves diagonally upward until it touches
the right leg above point $i$. From there it moves downward to point
$i+1$.  This peculiar ``law of refraction'' implies that only every
second point of the left half of the hypotenuse is connected to the
right half by a reflected path.  

We can apply the recursion again and write
\begin{displaymath}
  \begin{aligned}
  \lambda_1^{7} &= \quad\, \lambda_1^6 \quad + \quad \lambda_1^5 \\
                &= \lambda_1^5 + \lambda_1^3 + \lambda_1^4 + \lambda_1^1
  \end{aligned} 
\end{displaymath}
Here we have connected $6$ with $5$ along the hypotenuse and with $3$ along
a reflected path, and similarly for $5$. We iterate this path finding process
until all paths end on the left half of the hypotenuse (negative $t$). Here
the paths collect the initial values $\lambda_1^t=1$, hence 
$\lambda_1^7$ equals the number of different paths that connect the points
$-7,-5,\ldots,-1$ to the point $7$ on the hypotenuse. Instead of considering
each paths that starts on the left half of the hypotenuse separately we let
all paths start in the leftmost point $-7$. The rule for path finding then is:
if you are on an even index, move one unit to the right. If you are on an
odd index, there are two branches: one to the right, the other 45 degrees
upward and reflected down to the hypotenuse. Obviously, $\lambda_1^7$ equals the number of different paths that connects the leftmost point of $\Delta_8$
to the rightmost point according to this rules. Let $T_n(i)$ denote the
number of paths that connect the point $i$ with $n-1$ in $\Delta_n$.
Then we have
\begin{displaymath}
  T_n(-n+1) = \lambda_1^{n-1}\,.
\end{displaymath}
Now, starting at $-n+1$, we have two choices: move upward for a reflection that
will take us to point $1$ or move along the hypotenuse to point $-n+3$:
\begin{displaymath}
  T_n(-n+1) = T_n(1) + T_n(-n+3)\,.
\end{displaymath}
As we can see in Fig.~\ref{fig:proof} (top), the number of paths from $1$ to
$n-1$ is exactly the same as the total number of paths in $\Delta_{n/2}$.
Hence
\begin{displaymath}
  T_n(1) = T_{n/2}(-n/2+1)\,.
\end{displaymath}
Similarly, the number of paths from $-n+3$ to $n-1$ equals the
total number of paths in a slightly smaller triangle, as can be seen
in Fig.~\ref{fig:proof} (bottom). Hence we have
\begin{displaymath}
  T_n(-n+3) = T_{n-1}(-n+2)\,,
\end{displaymath}
and all three equations yield
\begin{displaymath}
   T_n(-n+1) = T_{n/2}(-n/2+1) + T_{n-1}(-n+2)\,.
\end{displaymath}
The derivation of a corresponding equation for odd values of $n$ is
straightforward. If we define
\begin{equation}
  \label{eq:def-F}
  F(n) := T_n(-n+1) = \lambda_1^{n-1}\,,
\end{equation}
the recursion for $T_n$ translates into the Fibonacci like recursion
\begin{equation}
  \label{eq:fibonacci}
  \begin{aligned}
  F(n) &= F(n-1) + F([n/2]) \\
  F(1) &= 1 
  \end{aligned}
\end{equation}
where $[x]$ refers to the integer part of $x$. 
The resulting sequence is known as
\htmladdnormallink{A033485}{http://www.research.att.com/~njas/sequences/A033485}
in \cite{oeis}. The generating function $g(z) = \sum_n F(n)\,z^n$
satisfies the functional equation
\begin{equation}
  \label{eq:generating-function}
  g(z)\,(1-z) = z + (1+z)\,g(z^2)\,,
\end{equation}
and is given by
\begin{equation}
  \label{eq:explicit-generating-function}
  g(z) = \frac{1}{2}\left(\frac{(1-z)^{-1}}{\prod_{k\geq 0}(1-z^{2^k})} - 1\right)\,.
\end{equation}
$F(n)$ can be evaluated numerically for values of $n$ that are larger
than the values feasible for simulations of LDM or for solving the rate
equation. The bottleneck for calculating $F(n)$ is memory, not CPU
time, since $n/2$ values must be stored to get $F(n)$. With 3
GByte of memory, we managed to calculate $F(n)$ for $n \leq
6\cdot10^8$. We will derive the asymptotics of $F(n)$ in the next section.
 
\begin{figure}
\includegraphics[angle=270,width=1.0\columnwidth]{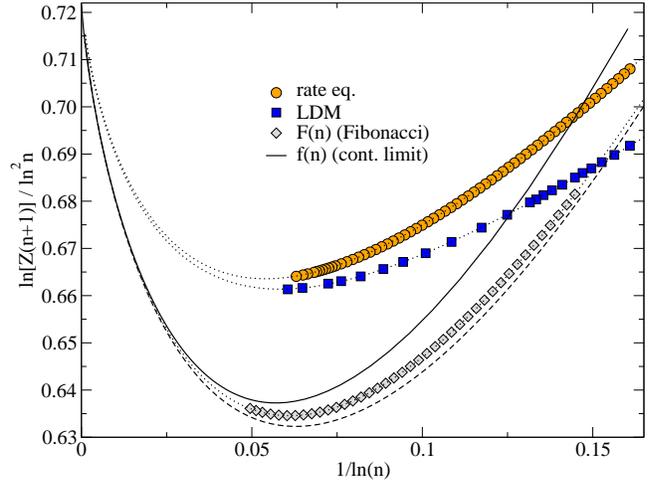}
\caption{Four models of LDM: Direct simulation ($Z=1/\ave{n\,L_n}$),
rate equation ($Z=\lambda_1^{n-1}$), the Fibonacci model $Z=F(n)$ from
\eqref{eq:fibonacci} and the similarity solution $Z=f(n)$ of the continuous rate
equation, given by \eqref{eq:def-littlef}. The dashed line represents \eqref{eq:fn-asympt}. All dotted lines are numerical fits of the type \eqref{eq:fit}.
\label{fig:allextrapolations}}
\end{figure}

Fig~\ref{fig:allextrapolations} shows $F(n)$ within the same scaling
as the simulations of LDM and the numerical solution of the rate
equation.  Apparently the similarity Ansatz does not capture the full
complexity of the LDM algorithm or the rate equation. Yet it
yields a very similar \emph{qualitative} behavior. And in the next section
we will show that
\begin{equation}
  \label{eq:Fn-limit}
  \lim_{n\to\infty}\frac{\ln F(n)}{\ln^2 n} = \frac{1}{2\ln 2}\,.
\end{equation}

\section{Continuum Limit}
\label{sec:continuum}

To analyze the rate
equations~(\ref{eq:lambda1eq}-\ref{eq:lambdainiteq}), it is convenient
to consider the continuum limit for $n\to\infty$. Asymptotically, a
continuum solution  may differ from the discrete problem in
corrections of order $1/n$. As we will see, such corrections are
inaccessible, as the asymptotic expansion is a series in terms of
$1/\ln(n)$.

\noindent
We rewrite Eq.~(\ref{eq:lambda1eq}) in terms of discrete differences,
\begin{equation}
\lambda_{i}^{t+1}-\lambda_{i}^{t}=-\left(\lambda_{i}^{t}-\lambda_{i-1}^{t}\right)
+\left(\lambda_{i}^{t}-\lambda_{i-1}^{t}+\lambda_{n-t}^{t}\right){\cal P}_{i}^t\,.
\label{eq:lambdabulk}
\end{equation}
Setting 
\begin{eqnarray}
t & = & sn\qquad(0\leq s\leq1)\,,\nonumber \\
i & = & xn\qquad(0\leq x\leq1-s)\,,
\label{eq:continuumlimit}\\
\lambda_{i}^{t} & = & y(x,s)\,,\nonumber 
\end{eqnarray}
we obtain for large $n$
\begin{eqnarray}
\frac{1}{n}\left[\frac{\partial}{\partial x}+\frac{\partial}{\partial s}\right]
y(x,s) & = & \Pi(x,s)\left[\frac{\partial}{n\partial x}y(x,s)+y(1-s,s)\right]
\,,
\label{eq:continuumbulk}
\end{eqnarray}
where we have set
\begin{equation}
{\cal P}_{i}^t\to\Pi(x,s)=\exp\left\{n\int_{0}^{x}d\xi\,\ln\alpha(\xi,s)\right\}
\label{eq:defPi}
\end{equation}
with
\begin{equation}
\alpha(x,s) = \frac{y(x,s)}{y(x,s)+y(1-s,s)}\leq1\,.
\label{eq:continuumalpha}
\end{equation}
The left-hand side of Eq.~(\ref{eq:continuumbulk}), as well as the
numerical solution of the full rate
equations~(\ref{eq:lambda1eq}-\ref{eq:lambdainiteq}) displayed in
Fig.~\ref{fig:lambdacontour}, again suggest a \emph{similarity Ansatz}
\begin{eqnarray}
y(x,s) & = & \gamma(s-x)\,.
\label{eq:continuumwaveAnsatz}
\end{eqnarray}
This Ansatz yields immediately for Eq.~(\ref{eq:continuumbulk}):
\begin{eqnarray}
0 & = & \Pi(x,s)\left[-\frac{1}{n}\,  \gamma\,'(s-x)+\gamma(2s-1)\right]\,.
\label{eq:remainderbulk}
\end{eqnarray}
For almost all $x>0$, the right-hand side vanishes by virtue of
$\Pi(x,s)\to0$, as indicated by Eq.~(\ref{eq:defPi}) for $\alpha<1$
and $n\to\infty$. Correspondingly, $\Pi(x=1-s,s)=0$ at the upper
boundary, which justifies the similarity solution for the continuum
limit of Eq.~(\ref{eq:lambda2eq}). Yet, $\Pi(x=0,s)=1$ for all $s$,
hence we are left with
\begin{equation}
\frac{1}{n}\, \gamma\,'(s) = \gamma(2s-1)\,,
\label{eq:continuumBClower}
\end{equation}
which can be interpreted as the continuous version of \eqref{eq:lambdafibo}.
From  the initial conditions of the discrete problem in
\eqref{eq:lambdainiteq} it is clear that $y(x,0)=1$. For the
similarity solution, this implies that 
\begin{equation}
\gamma(s)=1,\qquad(-1\leq s\leq0)\,.
\label{eq:continuumIC}
\end{equation}

Integrating \eqref{eq:continuumBClower}, we formally obtain
\begin{equation}
\gamma(s)  =  \gamma(0)+n\int_{0}^{s}d\xi\, \gamma(2\xi-1)\,.
\label{eq:inteq}
\end{equation}
Thus, we can evaluate the integral for $0\leq s\leq\frac{1}{2}$ to get
\begin{equation}
\gamma(s) =  1+ns,\qquad\left(0\leq s\leq\frac{1}{2}\right)\,.
\label{eq:f1}
\end{equation}
We can continue this process for $\frac{1}{2}\leq s\leq\frac{3}{4}$,
i.~e., $0\leq2s-1\leq\frac{1}{2}$, \emph{exactly} the domain of validity
of (\ref{eq:f1}), to obtain
\begin{eqnarray}
\gamma(s) & = & \gamma(0)+n\int_{0}^{\frac{1}{2}}d\xi\,\gamma(2\xi-1)+n\int_{\frac{1}{2}}^{s}d\xi\, 
\gamma(2\xi-1)\,,\nonumber \\
 & = & 1+ns+\frac{n^{2}}{4}\left(2s-1\right)^{2}
\qquad\left(\frac{1}{2}\leq s\leq\frac{3}{4}\right)\,.
\label{eq:f2}
\end{eqnarray}
The emergent pattern is best represented by defining
\begin{eqnarray}
\gamma_{k}(s) & = & \gamma(s),\qquad\left(1-2^{1-k}\leq s\leq1-2^{-k}\right)\,,
\label{eq:Deffk}
\end{eqnarray}
for $k=0,1,2,\ldots$, where Eqs.~(\ref{eq:continuumIC}-\ref{eq:f2})
represent $k=0$, $1$ and $2$. In general, we find that
\begin{eqnarray}
  \gamma_{k+1}(s) & = & \gamma_{k}\left(1-2^{-k}\right)+n\int_{1-2^{-k}}^{s}d\xi\,
  \gamma_{k}(2\xi-1)\,,
\label{eq:fkrecur}
\end{eqnarray}
which is solved by 
\begin{eqnarray}
\gamma_{k}(s) & = & \sum_{j=0}^{k}\frac{n^{j}}{j!\,2^{{j \choose
      2}}}\,\left(2^{j-1}s-2^{j-1}+1\right)^{j}\,.
\label{eq:fksolution}
\end{eqnarray}
For any $n$, we are interested in  $\gamma(s\to1)\sim\lim_{t\to n-1}\lambda_1^t$,
hence
\begin{equation}
\gamma(1) = \lim_{k\to\infty}\gamma_{k}\left(1-2^{-k}\right)=
\sum_{j=0}^{\infty}\frac{n^{j}}{j!\,2^{{j \choose 2}}}\,,
\label{eq:f1solution}
\end{equation}
which concludes our solution of \eqref{eq:continuumBClower}.
The sum for $\gamma(1)$ still depends on $n$, hence we define
\begin{equation}
  \label{eq:def-littlef}
  f(n) = \sum_{j=0}^{\infty}\frac{n^{j}}{j!\,2^{{j \choose 2}}}\,.
\end{equation}
Now $f(n)$ can be evaluated numerically for very large values of $n$.
Fig.~\ref{fig:allextrapolations} shows the result 
for $n \leq 2^{2000}$. Don't try this at home unless you have
a computer algebra system.
Interestingly, $\ln f(n)/\ln^2 n$ asymptotically approaches a value that is
extremely close to $1/2\ln 2$. In fact, an asymptotic analysis (see
Appendix) reveals
\begin{equation}
  \label{eq:fn-asympt}
  \begin{aligned}
  \frac{\ln \left[f(n)(n+1)\right]}{\ln^2 n} &\simeq \frac{1}{2\ln2} + \frac{1}{\ln n}\left(\frac{\ln\ln 2 + 1}{\ln{2}} +\frac{3}{2}\right) \\
   &+\frac{1}{\ln^2n}\left(\frac{\ln 2 + 4\ln\ln 2}{8} - \frac{\ln^2\ln 2}{2\ln2}\right) \\
   &-\frac{\ln\ln n}{\ln n}\,\frac{1}{\ln 2} - \frac{\ln\ln n}{\ln^2n} + 
    \frac{\ln^2\ln n}{\ln^2 n}\,\frac{1}{2\ln 2}\,,
  \end{aligned}
\end{equation}
which is the dashed line in Fig.~\ref{fig:allextrapolations}. The dotted
lines are numerical least square fits of the $\ln\ln n$ terms of this
scaling, i.e., fits of the form
\begin{equation}
  \label{eq:fit}
  \begin{aligned}
  \frac{\ln\left[Z(n)(n+1)\right]}{\ln^2 n} &\simeq \frac{1}{2\ln2} + \frac{1}{\ln n}\left(\frac{\ln\ln 2 + 1}{\ln{2}} +\frac{3}{2}\right) \\
   &+\frac{1}{\ln^2n}\left(\frac{\ln 2 + 4\ln\ln 2}{8} - \frac{\ln^2\ln 2}{2\ln2}\right) \\
   &+\frac{\ln\ln n}{\ln n}\,c_1 + \frac{\ln\ln n}{\ln^2n}\,c_2 + 
    \frac{\ln^2\ln n}{\ln^2 n}\,c_3\,.
  \end{aligned}
\end{equation}
with values for $c_i$ as shown shown in Table \ref{tab:fitparameter}.
\begin{table}
  \centering
  \begin{tabular}{c|cccc}
    $Z$ & $f$ & $F$ & $\lambda_1^{n-1}$ & $\ave{n\,L_n}^{-1}$ \\\hline
    $c_1$ & -1.44 & -1.45 & -1.22 & -1.24 \\
    $c_2$ & -1.00 & -1.42 & -3.06 & -3.86 \\
    $c_3$ & \phantom{-}0.72 & \phantom{-}1.01 & \phantom{-}1.23 & \phantom{-}1.55
  \end{tabular}
  \caption{Parameters for \eqref{eq:fit} used in Fig.~\ref{fig:allextrapolations}.}
  \label{tab:fitparameter}
\end{table}
Note that the series \eqref{eq:def-littlef} as a solution of
\eqref{eq:continuumBClower} and the first terms of the asymptotic
expansion \eqref{eq:fn-asympt} have been derived independently in the
context of dynamical systems \cite{moore:lakdawala:00}.

\section{Conclusion}
\label{sec:Conclusion}

The numerical data supports the claim that the complete statistics of
LDM is dominated by a single scale $\sim n^{-c\ln n}$, not just the
expectation as described in \eqref{eq:yakir}. The available data is
not sufficient to pin down the precise asymptotic scaling, however. In
fact a naive extrapolation of the available data even contradicts the
known asymptotic bound \eqref{eq:yakir2}. With its $\bigo{n\ln n}$
complexity, LDM is a very efficient algorithm, but probing the
asymptotics requires $\ln n$ to be large. This discrepancy of scales
eliminates simulations as a means to study the asymptotics of LDM and
calls for alternative approaches.

We have taken a step in the direction of a rigorous asymptotic
analysis by mapping the differencing algorithm onto a rate equation. 
The structure seen in the evolution of this rate equation (Fig.~\ref{fig:lambdacontour}) suggests a similarity Ansatz
\eqref{eq:similarity-discrete}. With the help of this Ansatz
we could reduce the exact recursion in
$\lambda$-space to the Fibonacci model \eqref{eq:fibonacci}.
The asymptotics of this model can be calculated, and it agrees with
\eqref{eq:yakir} and \eqref{eq:yakir2}. The same Ansatz plugged into the
rate equation even allows us to calculate the first
terms of an asymptotic expansion \eqref{eq:fn-asympt}. Although our
Ansatz does not yield a proof, the extracted asymptotic behavior
satisfies all previous constraints and provides a consistent
interpretation of the numerical results. Hence, our rate equations
pave the way for further systematic investigations.

{
\acknowledgement
We appreciate stimulating discussions George E. Hentschel and Cris Moore.
S.M. enjoyed the hospitality of the Cherry L. Emerson Center for
Scientific Computation at Emory University, where part of this work was done. Most simulations were run on the Linux-Cluster \textsc{Tina} at Magdeburg University.
S.M. was sponsored by the European Community's FP6 Information Society Technologies programme, contract IST-001935, EVERGROW.
}

\appendix

\section{Asymptotic Analysis}
\label{sec:proof-fibo}

To evaluate the series \eqref{eq:def-littlef} we apply Laplace's
saddle-point method for sums as described on p.~304 of Ref.~\cite{BO}.
For
\begin{displaymath}
a_{j} =  \frac{n^{j}}{j!\,2^{{j \choose 2}}}=e^{\phi_{j}}\,,
\end{displaymath}
the saddle point is determined by $D\phi_{j}=\phi_{j}-\phi_{j-1}=0$,
i.~e., $0=D\ln(a_{j})=\ln\left(a_{j}/a_{j-1}\right)$, or 
\begin{eqnarray}
1 & =\frac{a_{j}}{a_{j-1}} & =\frac{n}{j\,2^{j-1}}\,.
\label{eq:saddlepointcondition}
\end{eqnarray}
Hence, we obtain a moving ($n$-dependent) saddle point at
\begin{eqnarray}
j_{0} & \sim &
\frac{\ln n}{\ln 2}-\frac{\ln\left(\frac{\ln n}{\ln 2}\right)}{\ln 2}+1+\frac{\ln\left(\frac{\ln n}{\ln 2}\right)}{\ln 2\ln n}-\frac{1}{\ln n}+\ldots,
\label{eq:saddlepoint}
\end{eqnarray}
including terms to the order needed to determine $f(n)$
up to the correct prefactor. We keep
the $1/\ln(n)$-corrections, since $\phi_{j}$ contains terms like
$j_{0}\ln(n)$. In particular, it is
\begin{eqnarray}
\phi_{j} & = & j\,\ln n -\ln j!-\frac{j(j-1)}{2}\ln 2\,.
\label{eq:phi}
\end{eqnarray}
As the saddle point $j_{0}$ is large for large $n$, we can replace
$j!$ by its Stirling-series~\cite{BO}.
Then, we expand around the saddle point by substituting $j=j_{0}+\eta$,
keeping only terms to \emph{2nd} order in $\eta$ and those that are
non-vanishing for $n\to\infty$. We find
\begin{displaymath}
\phi_{j_{0}+\eta} \sim \frac{\ln^2 n}{2\ln 2}-\frac{1}{2}\ln(2\pi)-\frac{\ln 2}{2}\eta(\eta+1) + \mathcal{C}(n)\,,
\end{displaymath}
with log-polynomial corrections
\begin{equation}
\begin{aligned}
{\cal C}(n) \sim &
-\frac{\ln n}{\ln 2}
\left[\ln\left(\frac{\ln n}{\ln 2}\right)-1-\frac{\ln 2}{2}\right]\\
&+\left[\frac{1}{2\ln 2}\ln^2\left(\frac{\ln n}{\ln(2)}\right)-\ln\left(\frac{\ln n}{\ln(2)}\right)\right]\,.
\end{aligned}
\label{eq:logpoly}
\end{equation}
We finally obtain for the asymptotic expansion of (\ref{eq:fn-asympt}):
\begin{eqnarray}
f(n) & \sim & \int_{-\infty}^{\infty}d\eta\,\exp\left(\phi_{j_{0}+\eta}\right) \nonumber\\
 & \sim & \frac{2^{\frac{1}{8}}}{\sqrt{\ln2}}\,\exp\left\{
\frac{\ln^2 n}{2\ln2}+{\cal C}(n)\right\}
\label{eq:f1asympt}
\end{eqnarray}

\bibliographystyle{unsrt}
\bibliography{complexity,cs,math,mertens}

\end{document}